# The Dyer-Lashof Algebra in Bordism

Terrence Bisson    André Joyal
bisson@canisius.edu    joyal@math.uqam.ca

We present a theory of Dyer-Lashof operations in unoriented bordism (the canonical splitting $N_*(X) \simeq N_* \otimes H_*(X)$, where $N_*( )$ is unoriented bordism and $H_*( )$ is homology mod 2, does not respect these operations). For any finite covering space we define a "polynomial functor" from the category of topological spaces to itself. If the covering space is a closed manifold we obtain an operation defined on the bordism of any $E_\infty$-space. A certain sequence of operations called squaring operations are defined from two-fold coverings; they satisfy the Cartan formula and also a generalization of the Adem relations that is formulated by using Lubin's theory of isogenies of formal group laws. We call a ring equipped with such a sequence of squaring operations a *D-ring*, and observe that the bordism ring of any free $E_\infty$-space is free as a *D*-ring. In particular, the bordism ring of finite covering manifolds is the free *D*-ring on one generator. In a second compte-rendu we discuss the (Nishida) relations between the Landweber-Novikov and the Dyer-Lashof operations, and show how to represent the Dyer-Lashof operations in terms of their actions on the characteristic numbers of manifolds.

## 1. The algebra of covering manifolds.

We begin with the observation that a covering space $p : T \to B$ can be used to define a functor $X \mapsto p(X)$ from the category of topological spaces to itself, where

$$p(X) = \{(u, b) \mid b \in B,\ u : p^{-1}(b) \to X\}.$$

Then $p(X)$ is the total space of a bundle over $B$ with fibers $X^{p^{-1}(b)}$, and any continuous map $f : X \to Y$ induces a continuous map $p(f) : p(X) \to p(Y)$. We shall say that $p( )$ is a *polynomial functor*. For functors $F$ and $G$ from the category of topological spaces to itself, we have functors $F + G$, $F \times G$ and $F \circ G$ given by $(F + G)(X) = F(X) + G(X)$, $(F \times G)(X) = F(X) \times G(X)$, and $(F \circ G)(X) = F(G(X))$. Polynomial functors happen to be closed under these operations, and we obtain well-defined operations $p + q$, $p \times q$ and $p \circ q$ on coverings. These operations satisfy the kinds of identities that one should expect for an algebra of polynomials.

We define the *derivative* $p'$ of a covering $p : T \to B$ to be the covering whose base space is $T$ and whose fiber over $t \in T$ is the set $p^{-1}(p(t)) - \{t\}$. The rules of differential calculus apply: $(p + q)' = p' + q'$, $(p \times q)' = p' \times q + p \times q'$ and $(p \circ q)' = (p' \circ q) \times q'$. If we observe that the total space of $p$ is $p'(1)$ (where 1 denotes a single point) and that its base space is $p(1)$ the formula $(p \times q)'(1) = p'(1) \times q(1) + p(1) \times q'(1)$ expresses the total space of $(p \times q)$ in terms of the total and based spaces of $p$ and $q$. Similarly for the formula $(p \circ q)'(1) = p'(q(1)) \times q'(1)$.

**Remark 1:** There is a parallel between this algebra of covering spaces and the algebra of combinatorial species developed in [9] and [10].

**Remark 2:** By using the Euler-Poincare characteristic one can associate a polynomial $\chi(p)$ to any covering $p$ of a finite complex. We have $\chi(p+q) = \chi(p) + \chi(q)$, $\chi(p \times q) = \chi(p) \times \chi(q)$, $\chi(p \circ q) = \chi(p) \circ \chi(q)$, and $\chi(p') = \chi(p)'$.

**Remark 3:** It is also possible to define various kinds of higher differential operators on coverings. For example, the group $\Sigma_2$ acts on any second derivative $p''$ by permuting the order of differentiation, and we can define

$$\frac{1}{2!} \frac{d^2 p}{dx^2} = p''/_{\Sigma_2}.$$

Higher divided derivatives can be handled similarly.

**Remark 4:** Polynomial functors of $n$ variables are easily defined. They are obtained from $n$-tuples $(p_1, \ldots, p_n)$ where $p_i : T_i \to B$ is a finite covering for every $i$.

Let us now consider coverings of smooth compact manifolds. We say that two coverings of closed manifolds are *cobordant* if together they form the boundary of a covering. Let $N_*\Sigma$ denote the set of cobordism classes of closed coverings. Let $N_d\Sigma_n$ denote the set of cobordism classes of degree $n$ (i.e. $n$-fold) coverings over closed manifolds of dimension $d$.

**Proposition 1.** *The operations of sum $+$, product $\times$, and composition $\circ$ are compatible with the cobordism relation on closed coverings. They define on $N_*\Sigma$ the structure of a commutative $\mathbf{Z}_2$ algebra, graded by dimension.*

Notice that if $p \in N_k\Sigma_m$ and $q \in N_r\Sigma_n$ then $p \circ q \in N_{mr+k}\Sigma_{mn}$. This defines in particular an action of $N_*\Sigma$ on $N_*\Sigma_0 = N_*$. More generally, let us see that $N_*\Sigma$ acts on the bordism ring of any $E_\infty$-space.

Recall (see [1], [18]) that an $E_\infty$-space $X$ has structure maps $E\Sigma_n \times_{\Sigma_n} X^n \to X$ for each $n$. These structure maps give rise to structure maps $p(X) \to X$ for every degree $n$ covering space $p : T \to B$. To see this it suffices to express $p$ as a pull back of the tautological $n$-fold covering $u_n$ of $B\Sigma_n$ along some map $B \to B\Sigma_n$. This furnishes a map $p(X) \to u_n(X) = E\Sigma_n \times_{\Sigma_n} X^n$ and the structure map $p(X) \to X$ is then obtained by composing with $u_n(X) \to X$.

Recall (see [6] for instance) that an element of $N_*X$ is the bordism class of a pair $(M, f)$ where $f : M \to X$ and $M$ is a compact manifold; then $p(M)$ is a compact manifold and the structure map for $X$ gives $p(M) \to p(X) \to X$, representing an element in $N_*X$.

**Proposition 2.** *Let $X$ be an $E_\infty$-space. Each covering of degree $n$ and dimension $d$ defines an operation $N_m X \to N_{nm+d} X$. Cobordant covering spaces give the same operation. Moreover, for double coverings these operations are additive.*

It should be noted that tom Dieck [7] and Alliston [3] develop bordism Dyer-Lashof operations which agree with ours; the relationship will be clearer after section 2.

**Example:** The classifying space for finite coverings is $B\Sigma_*$ the disjoint union of the classifying spaces of the symmetric groups $B\Sigma_n$. Then $N_*B\Sigma_* = N_*\Sigma$ and $B\Sigma_*$ has a natural $E_\infty$-space structure defined from disjoint sum. The covering operations on $N_*B\Sigma_*$ correspond to composition of coverings.

**Remark:** It is a classical result [19], [8], [12] that the inclusion $i : \Sigma_{n-1} \subset \Sigma_n$ defines a *split monomorphism* $i_* : N_*\Sigma_{n-1} \to N_*\Sigma_n$. In our setting $i_*$ is the map $p \mapsto x \times p$. It is easy to see, by applying the rules of differential calculus, that the map

$$q \mapsto \frac{dq}{dx} + x\frac{1}{2!}\frac{d^2q}{dx^2} + x^2\frac{1}{3!}\frac{d^3q}{dx^3} + \cdots$$

is a splitting [11].

For any space $X$ let $\epsilon : N_*(X) \to H_*(X)$ denote the Thom reduction, where $H_*(\ )$ is mod 2 homology. If $(M, f) \in N_*(X)$ we have $\epsilon(M, f) = f_*(\mu_M)$ where $\mu_M$ denotes the fundamental homology class of $M$. If $X$ is an $E_\infty$-space then each covering of degree $n$ and dimension $d$ defines an operation $H_m X \to H_{nm+d} X$ which is the Thom reduction of the corresponding operation in bordism.

We now describe the sequence of cobordism class of double coverings that leads to the concept of $D$-rings. It is a classical result that $N^*(RP^\infty) = N_*[[t]]$. Let $q_k$ in $N_*B\Sigma_2 = N_*(RP^\infty)$ be represented by the canonical inclusion $RP^k \hookrightarrow RP^\infty$. The sequence $q_0, q_1, \ldots$ is a basis of the $N_*$-module $N_*(RP^\infty)$. The Kronecker pairing $N^*(RP^\infty) \times N_*(RP^\infty) \to N_*$ defines an exact duality between $N^*(RP^\infty)$ and $N_*(RP^\infty)$. Let $d_0, d_1, \ldots$ be the basis dual to the basis $t^0, t^1, t^2, \ldots$ under the Kronecker pairing. The relation between the two bases of $N_*(RP^\infty)$ can be expressed as an equality of generating series

$$(\sum_{i \geq 0}[RP^i]x^i)(\sum_{k \geq 0} d_k x^k) = (\sum_{n \geq 0} q_n x^n),$$

where $x$ is a formal indeterminate. We have $d_0 = q_0$, and $d_1 = q_1$ since $[RP^0] = 1$ and $[RP^1] = 0$. It turns out (see [2] for instance) that $d_n$ can be represented by the Milnor hypersurface $H(n,1) \hookrightarrow RP^n \times RP^1 \to RP^n$. The coverings $d_n$ and $q_n$ give operations which are distinct in bordism but agree in mod 2 homology.

## 2. $D$-rings and Dyer-Lashof operations

Recall that a formal group law over a commutative ring $R$ is a formal power series $F(x, y) \in R[[x, y]]$ which satisfies identities corresponding to associativity and unit; (see Quillen [21] or Lazard [13] for instance). We say that a formal group law $F$ has order two if $F(x, x) = 0$.

The Lazard ring (for formal group laws of order two) is the commutative ring with generators $a_{i,j}$ and relations making $F(x, y) = \sum a_{i,j} x^i y^j$ a formal group law of order two. Let us temporarily denote this Lazard ring by $L$. Then for any ring $R$ and any formal group law $G(x, y) \in R[[x, y]]$ of order two there is a unique ring homomorphism $\phi : L \to R$ such that $(\phi F)(x, y) = G(x, y)$. Quillen [21] showed that $L$ is naturally isomorphic to $N_* = N_*(pt)$. This provides a beautiful interpretation of Thom's original calculation of the unoriented cobordism ring.

Let $R$ be a commutative ring and let $F \in R[[x, y]]$ be a formal group law of order two (this implies that $R$ is a $\mathbf{Z}_2$-algebra). According to Lubin [14] there exists a unique formal group law $F_t$ defined over $R[[t]]$ such that $h_t(x) = xF(x, t)$ is a morphism $h_t : F \to F_t$. The kernel of $h_t$ is $\{0, t\}$, which is a group under the $F$-addition $x +_F y = F(x, y)$. We will refer to $F_t$ as the *Lubin quotient* of $F$ by $\{0, t\}$ and to $h_t$ as the *isogeny*. The construction can be iterated and a Lubin quotient $F_{t,s}$ of $F_t$ can be obtained by further killing $h_t(s) \in R[[t, s]]$. The composite isogeny $F \to F_t \to F_{t,s}$ is

$$h_{t,s}(x) = h_t(x) F_t(h_t(x), h_t(s)) = xF(x, t)F(x, s)F(x, F(s, t))$$

Its kernel consists of $\{0, t, s, F(s, t)\}$, which is an elementary abelian 2-group under the $F$-addition. By doing the construction in a different order we obtain $F_{s,t}$ but it turns out that $F_{t,s} = F_{s,t}$.

**Definition:** A *D-ring* is a commutative ring $R$ together with a formal group law of order two $F$ defined over $R$ and a ring homormorphism $D_t : R \to R[[t]]$ called the *total square*, satisfying the following conditions:
  i) $D_0(a) = a^2$ for every $a$ in $R$;
  ii) $D_t(F) = F_t$;
  iii) $D_t \circ D_s$ is symmetric in $t$ and $s$. Here we have extended $D_t : R \to R[[t]]$ to $D_t : R[[s]] \to R[[s, t]]$ by defining $D_t(s) = h_t(s) = sF(s, t)$.

A *morphism* of $D$-rings is a ring homomorphism which preserves the formal group laws and the total squares. A $D$-ring is also an algebra over the Lazard ring $N_*$, and a morphism of $D$-rings is a morphism of $N_*$-algebras.

A $D$-ring is *graded* if $R$ is graded and $F$ is homogeneous in grade $-1$ and $D_t(x)$ has grading $2i$ in $R[[t]]$ for each element of grading $i$ in $R$ (where $t$ and $s$ have grading $-1$).

**Example:** The Lazard ring $N_*$ has a unique ring homomorphism $D_t : N_* \to N_*[[t]]$ such that $D_t(F) = F_t$, and this defines a $D$-structure on $N_*$. Thus $N_*$ is initial in the category of $D$-rings.

**Proposition.** *If $X$ is an $E_\infty$-space then $N_*X$ is a commutative ring under Pontryagin product; it is also an $N_*$-algebra. If $d_0, d_1, \ldots$ are the double coverings described in the previous section then the total squaring*

$$D_t(x) = \sum_n d_n(x) t^n$$

*gives an $D$-structure on $N_*X$.*

**Example:** $BO_*$, the disjoint union of the classifying spaces of the orthogonal groups $BO(n)$, is an $E_\infty$-space with $N_*BO_* = N_*[b_0, b_1, \ldots]$. It forms a $D$-ring with $F$ given by the cobordism formal group law over $N_*$ and with $D_t$ determined by

$$D_t(b)(xF(x, t)) = b(x)b(F(x, t))$$

where $b(x) = \sum b_i x^i$.

We shall refer to any $D$-ring $R$ with $F = (+)$ as a *Q-ring*. The mod 2 homology of an $E_\infty$-space $E$ is a $Q$-ring, and the Thom reduction $\epsilon : N_*(E) \to H_*(E)$ is a morphism of $D$-rings.

**Proposition.** *A commutative ring $R$ is a $Q$-ring if and only if it has a sequence of additive operations $q_n : R \to R$ which satisfy the following three conditions:*
  i) *Squaring: $q_0(x) = x^2$ for all $x \in R$.*
  ii) *Cartan formula: $q_n(xy) = \sum_{i+j=n} q_i(x) q_j(y)$ for all $x, y \in R$.*
  iii) *Adem relations: $q_m(q_n(x)) = \sum_i \binom{i-n-1}{2i-m-n} q_{m+2n-2i}(q_i(x))$ for all $x \in R$.*
*In the graded case, $\mathrm{grade}(q_n(x)) = 2 \cdot \mathrm{grade}(x) + n$.*

This is exactly an action of the classical Dyer-Lashof algebra on $R$. This idea of writing Adem relations via generating series is suggested by [4] and by Bullett and MacDonald [5]. See [17], [15], [16] for background on Dyer-Lashof operations.

**Example:** The $Q$-structure on $H_* BO_* = \mathbf{Z}_2[b_0, b_1, \ldots]$ is characterized by

$$Q_t(b)(x(x+t)) = b(x) b(x+t)$$

where $b(x) = \sum b_i x^i$. This expresses via generating series a calculation of Priddy's in [20].

Notice that if $A$ and $B$ are $Q$-rings then $A \otimes_{N_*} B = A \otimes_{\mathbf{Z}_2} B = A \otimes B$ is a $Q$-ring. Let $Q\langle M \rangle$ denote the free $Q$-ring generated by a $\mathbf{Z}_2$-vector space $M$. If $M$ has a comultiplication, then $Q\langle M \rangle$ has a comultiplication extending it which is a morphism of $Q$-rings.

Recall that $E_\infty(X)$ is the free $E_\infty$-space generated by $X$ (see [18] or [1] for background). The following is a classical result:

**Theorem 1.** (May [17]) *For any space $X$ the canonical map*

$$Q\langle H_* X \rangle \to H_* E_\infty(X)$$

*is an isomorphism which preserves the comultiplication. In particular, $H_* B\Sigma_* = Q\langle x \rangle$ is the free $Q$-ring on one generator.*

If $A$ and $B$ are $D$-rings then $A \otimes_{N_*} B$ is naturally a $D$-ring. Let us denote $D\langle M \rangle$ denote the $D$-ring freely generated by an $N_*$-module $M$. If $M$ is a coalgebra in the category of $N_*$-modules, then $D\langle M \rangle$ has a comultiplication.

**Theorem 2.** *The bordism of an $E_\infty$-space is an $D$-ring. Moreover, for any space $X$ the canonical map*

$$D\langle N_* X \rangle \to N_* E_\infty(X)$$

*is an isomorphism which preserves the comultiplication. In particular, $N_* \Sigma = N_*(B\Sigma) = D\langle x \rangle$ is the free $D$-ring on one generator.*

Thus, both $D\langle x \rangle$ and $N_* \Sigma$ are algebras equipped with operations of substitution; the former because it is the set of unary operations in the theory of $D$-rings and the latter because we have defined a substitution operation among coverings of manifolds. The above theorem says that the canonical isomorphism of $D$-rings $D\langle x \rangle \to N_* \Sigma$ which sends the generator $x$ to the unique non-zero element $x$ in $N_0(B\Sigma_1)$ preserves the operations of substitution.

## References

[1] J. F. Adams. Infinite loop spaces. Ann. of Math. Studies no. 90, Princeton 1978.

[2] M. A. Aguilar, Generators for the bordism algebra of immersions, Trans. Amer. Math. Soc. 316 (1989).

[3] R. M. Alliston, Dyer-Lashof operations and bordism. Ph.D. Thesis, University of Virginia, 1976.

[4] T. P. Bisson, Divided sequences and bialgebras of homology operations. Ph.D. Thesis, Duke 1977.

[5] S. R. Bullett, I. G. Macdonald, On the Adem relations, Topology 21 (1982), 329-332.

[6] P. E. Conner (and E. E. Floyd). Differentiable Periodic Maps (2nd ed.), Lect. Notes in Math. no. 738, Springer 1979.

[7] T. tom Dieck, Steenrod operationen in Kobordismen-Theorie, Math. Z. 107 (1968), 380-401.

[8] A. Dold, Decomposition Theorem for $S(n)$-complexes, Ann. of Math.(2) 75 (1962), 8-16.

[9] A. Joyal, Une théorie combinatoire des séries formelles, Adv. in Math. 42 (1981), 1-82.

[10] A. Joyal, Foncteurs analytiques et espèces de structures, Lect. Notes in Math. no. 1234, Springer 1985.

[11] A. Joyal, Calcul Integral Combinatoire et homologie des groupes symetriques, C. R. Math. Rep. Acad. Sci. Canada 12 (1985).

[12] D. S. Kahn, S. B. Priddy, Application of transfer to stable homotopy theory, Bull. AMS 78 (1972), 981-987.

[13] M. Lazard. Commutative Formal Groups, Lect. Notes in Math. no. 443, Springer 1975.

[14] J. Lubin, Finite subgroups and isogenies of one-parameter formal Lie groups, Ann. of Math. 85 (1967) 296-302.

[15] I. Madsen, On the action of the Dyer-Lashof algebra in $H_*(G)$, Pacific J. Math. 60 (1975), 235-275.

[16] I. Madsen, R. J. Milgram. The classifying spaces for surgery and cobordism of manifolds, Ann. of Math. Studies no. 92, Princeton 1979.

[17] J. P. May, Homology operations on infinite loop spaces, in Proc. Symp. Pure Math. 22, Amer. Math. Soc. 1971, 171-185.

[18] J. P. May, Infinite loop space theory, Bull. Amer. Math. Soc. 83 (1977), 456-494.

[19] M. Nakaoka, Decomposition theorems for the homology of the symmetric groups, Ann. of Math. 71 (1960), 16-42.

[20] S. B. Priddy, Dyer-Lashof operations for the classifying spaces of certain matrix groups, Quart. J. Math. Oxford 26 (1975), 179-193.

[21] D. Quillen, Elementary proofs of some results of cobordism theory using Steenrod operations, Adv. in Math. 7 (1971), 29-56.

(*) Canisius College, Buffalo, N.Y. (U.S.A). e-mail: *bisson@canisius.edu*.

(**) Département de Mathématiques, Université du Québec à Montréal, Montréal, Québec H3C 3P8. e-mail: *joyal@math.uqam.ca*.